\documentclass[12pt,leqno]{article}
\usepackage{amsmath,amsfonts,amscd,ProtocoloTope,url}

\title{\Large\bf   Tautness  for riemannian foliations \\ on non-compact manifolds}

\author{
JosŽ Ignacio Royo Prieto\thanks{Departamento de Matem‡tica Aplicada.
Escuela Superior de Ingenier'a.  Universidad del Pa's Vasco.  Alameda
de Urquijo s/n.  48013 Bilbao.
Espa–a.  {\sl joseignacio.royo@ehu.es}.  Partially supported by the UPV - EHU
grant 127.310-E-14790/2002 by a PostGrant from the Gobierno Vasco -
Eusko Jaurlaritza and by the MCyT of the Spanish Government.} \\
{\small Universidad del Pa's Vasco.} \and Martintxo
Saralegi-Aranguren\thanks{UPRES-EA 2462 LML. FacultŽ Jean Perrin. 
UniversitŽ d'Artois.  Rue Jean Souvraz SP 18.  62 307 Lens Cedex -
France.  {\sl saralegi@euler.univ-artois.fr}.  Partially supported by
the UPV - EHU grant 127.310-E-14790/2002 } \\
{\small UniversitŽ d'Artois } \and Robert Wolak\thanks{Instytut
Matematyki.  Uniwersytet Jagiellonski.  Wl.  Reymonta4, 30 059 Krakow
- Poland.  {\sl wolak@im.uj.edu.pl} Partially supported by the KBN
grant 2 PO3A 021 25.} \\ {\small Uniwersytet Jagiellonski} }

\begin{document}  
\maketitle

\begin{abstract}
 For a riemannian foliation $\F$ on a closed manifold $M$, it is known
    that $\F$ is taut (i.e. the leaves are minimal submanifolds) if and
    only if the (tautness) class  defined by 
    the mean curvature form $\kappa_\mu$ (relatively to a suitable riemannian 
    metric $\mu$) is
    zero  (cf. \cite{Suso}). In the transversally orientable case, tautness is
    equivalent to the non-vanishing of the top basic cohomology group
    $\hiru{H}{n}{\mf}$, where $n = \codim \F$   (cf. \cite{MA}). By the PoincarŽ Duality 
   (cf. \cite{KT2}) this last condition is equivalent to the non-vanishing of the 
     basic twisted cohomology group
    $\lau{H}{0}{\kappa_\mu}{\mf}$, when $M$ is oriented.

    When $M$ is not compact, the tautness class is not even defined in general. In
    this work, we recover the previous study and results for a 
    particular case of riemannian
    foliations on non compact manifolds: the regular
    part of a singular riemannian foliation on a compact manifold
    (CERF).
    \end{abstract}

    The study of taut foliations (foliations where all the leaves are minimal 
    submanifolds for some riemannian metric) 
    has been an important part of the research in
    regular foliations on riemannian manifolds:  F. Kamber and Ph. 
    Tondeur (cf. \cite{KT}), H. Rummler (cf. \cite{RU}), as well as  D. Sullivan (cf. 
    \cite{SU})  were 
    the first  ones to present algebraical or variational characterizations of such
    foliations. A. Haefliger's paper \cite{HA} proved to be  an important step in
    the development of the theory. He showed that ``being taut" is a transverse
    property, i.e. it depends only on the holonomy pseudogroup of the 
    foliation $\F$.
    This led Y. Carrire (cf. \cite{CA}) to propose a characterization of taut riemannian
    foliations on a compact manifold $M$ as those foliations for which the top
    dimensional basic cohomology group $\hiru{H}{n}{\mf}$ is non-trivial, 
    i.e. isomorphic to 
    $\R$, being $n = \codim \F$.  For  a concise presentation of the history 
    of the basic cohomology and tautness we
    refer to V. Sergiescu's appendix \cite{SER_A} in \cite{Mo},
    which shows a close relation between the finiteness of basic cohomology, 
    PoincarŽ Duality
    Property in basic cohomology and tautness. 
    
    The problem was positively solved by X. Masa in \cite{MA}, in the 
    following way: if $M$ is compact and oriented  and $\F$ is 
    riemannian and transversally oriented, then 
    \begin{equation}
	\label{Masa}
	\F \hbox{ is taut if and 
    only if } \hiru{H}{n}{\mf} = \R.
    \end{equation}
    
    The work of A. çlvarez (cf. \cite{Suso}) removes the 
    orientability condition on $M$ and 
    gives another 
    characterization of taut riemannian
    foliations on a compact manifold. He constructs a cohomological 
    class $\kappab \in \hiru{H}{1}{\mf}$, the tautness class, whose 
    vanishing is equivalent to the tautness of $\F$. This 
    class is defined from  the mean curvature form $\kappa_\mu$ of a
    bundle-like metric $\mu$, which we can suppose to be basic due to D.
    Dom'nguez (cf.  \cite{D}).
    
    We have a third cohomological characterization of the tautness 
    of $\F$. The PoincarŽ Duality of \cite{KT2} implies that this 
    property is equivalent to the non-vanishing of  the  basic twisted cohomology group
    $\lau{H}{0}{\kappa}{\mf}$ when $M$ is orientable.
    
    \medskip
    
    The situation is more complicated when the manifold $M$ 
    is not compact. For example, the mean curvature form $\kappa_{\mu}$ may 
    be a basic form without being closed (cf. \cite{CE}).
    
    \smallskip
    
    We consider in this work a particular case of a riemannian 
    foliation $\F$ on a non-compact manifold $M$: 
    the Compactly Embeddable Riemannian Foliations
    or CERFs.
In this context,
we have a compact manifold $N$ endowed with a singular riemannian 
foliation $\HH$ (in the sense of \cite{Mo}) in such a way that $M$ is the regular stratum of 
$N$ with $\F = 
\HH|_{M}$.
We consider a class of bundle-like metrics on $M$  for which we construct a
tautness class $\kappab = [\kappa_{\mu}] \in \hiru{H}{1}{\mf}$ which
is independent of the choice of $\mu$.  We prove that the tautness of
$\F$ is equivalent to any of the following three properties:
    \begin{itemize}
	\item[-]     $\kappab = 0$, 
	\item[-] $\lau{H}{0}{\kappa_{\mu}}{\mf} \ne 0$,
	\item[-]     $\lau{H}{n}{c}{\mf}= \R$, where $n = \codim \F$, when 
	$\F$ is transversally oriented.

	\end{itemize}
 Notice that in the second characterization we have eliminated the 
 orientation hypothesis.
    We also prove that the cohomology groups 
    $\lau{H}{0}{\kappa}{\mf} $ and $\lau{H}{n}{c}{\mf}$ are $0$ or 
    $\R$. 
    
The standard method to prove the equivalence of the second and third 
conditions is to use the PoincarŽ Duality Property (PDP in short). 
However, although reasonable, the PDP for the basic cohomology has been
proved neither for $(M,\F)$ nor for $(N,\HH)$.
So, we shall proceed by proving that the first condition is 
equivalent to
 either 
of the two remaining ones. 
The proof of the equivalence of the first and second condition is purely 
algebraic. To obtain the second equivalence we use Molino's 
desingularisation $(\wt{N},\wt\HH)$ of $(N,\HH)$.  The key point is the
following.
    $$
   \F \hbox{ is taut } \Longleftrightarrow \wt\HH \hbox{ is taut. }
   $$
  Note that this equivalence cannot  be extended to the singular riemannian 
  foliation
  $\HH$ itself,   since the existence of 
    leaves with different dimensions implies the non-existence of ``minimal 
    metrics." Then the comparison of the corresponding basic 
    cohomology groups completes the proof. In this way we have avoided 
    in the proof any reference to  the PDP of the  basic cohomology of 
    the foliated 
    manifolds we are studying.

		\bigskip

	In the sequel $M$ and $N$  are connected, second countable, Haussdorff,   
	without boundary and smooth 
	(of class $C^\infty$) manifolds  of dimension	$m$. 
	All the maps considered are smooth unless something else is indicated.

	\smallskip

     \section{Riemannian 
     foliations\protect\footnote{For the notions related 
  to riemannian foliations we refer the reader to \cite{Mo,To}.}}
    
    The framework category of this work is that of CERFs. They 
    are riemannian foliations embedded in singular riemannian 
    foliations on compact manifolds. Before introducing  this notion, we need to 
    recall some important facts about singular riemannian foliations.
    
    \prg{\bf The SRFs }. A {\em{\rm S}ingular {\rm R}iemannian 
    {\rm F}oliation\protect\footnote{For the notions related 
   to singular riemannian foliations we refer the reader to
   \cite{BM,Mo0,Mo,Mo1}.}} (SRF for short) on a manifold $N$ is a
   partition $\mathcal{H}$ by connected immersed submanifolds, called
   {\em leaves}, verifying the following properties:

    \begin{itemize}
    \item[I-] The module of smooth vector fields tangent to the leaves is 
    transitive on each leaf.

    \item[II-] There exists a riemannian metric $\mu$ on $N$, called {\em 
    adapted metric},  such that each geodesic 
    that is perpendicular at one point to a leaf remains perpendicular to every 
    leaf it meets.
    \end{itemize}
    The first condition implies that $(N,\mathcal{H})$ is a singular foliation 
    in the sense of \cite{St} and \cite{Su}. Notice that the 
    restriction of $\mathcal{H}$ to a saturated open subset 
    produces an SRF. Each (regular) {\rm R}iemannian {\rm F}oliation (RF for
    short) is an SRF, but the first interesting examples are the following:
    \begin{itemize}
    \item[-] The orbits of the action by isometries of a Lie group.
    \item[-] The closures of the leaves of a regular riemannian foliation.
    \end{itemize}
    
    \prg{\bf Stratification}.  Classifying the points of
  $N$ following the dimension of the leaves one gets a stratification
  $\SH$ of 
  $N$ whose elements are called {\em strata}. The restriction of 
  $\mathcal{H}$ to a stratum $S$ is a RF
  $\mathcal{H}_{S}$. The strata are ordered by: $S_{1} \preceq S_{2} 
  \Leftrightarrow S_{1}\subset \overline{S_{2}}$. The minimal (resp. 
  maximal) strata are the closed strata (resp. open strata). Since $N$ is 
  connected, there is just one open stratum, denoted 
  $R_{\mathcal{H}}$. It is a dense subset.
  This is the {\em regular stratum}, the other strata are the {\em 
  singular strata}. The family of singular strata is written 
  $\SH^{^{sin}}$.
  The {\em dimension of the foliation} $\mathcal{H}$ is the dimension of the 
  biggest leaves of $\mathcal{H}$.
  
The  {\em depth }  of $\SH$, written $ \depth  \SH$, is defined to be the largest 
 $i$ for which there 
 exists a chain of strata $S_0  \prec S_1 \prec \cdots \prec S_i$. So, 
 $ \depth  \SH= 0$ if and only if the foliation 
 ${\mathcal  H}$ is regular. 

 The  {\em depth}  of  a 
 stratum $S \in \SH$ , written $ \depth_{\HH} S$, is defined to be the largest 
  $i$ for which there 
  exists a chain of strata $S_0  \prec S_1 \prec \cdots \prec S_i= S$. So, 
  $ \depth_{\HH}  S= 0$ (resp. $\depth_{\HH}  S=\depth \SH$) if and only if the stratum $S$ is 
  minimal (resp. regular).
 For each $i  \in \Z$ we write
 $$
 \Sigma_{i} = \Sigma_{i}(N) = \cup \{ S \in \SH \ | \ \depth _{\HH} S \leq 
 i\}.
 $$
We have $\Sigma_{<0}=\emptyset$, $\Sigma_{\depth \Sh -1} = N\menos 
R_{\HH}$ and $\Sigma_{i} = N$ if $i\geq\depth \SH$.  
The union of closed (minimal) strata is  $ \Sigma_0$.

    \prg {\bf The CERFs.}  Consider a riemannian foliation $\F$ on a 
    manifold $M$. We say that $\F$ 
    is a {\em {\rm C}ompactly 
    {\rm E}mbeddable {\rm R}iemannian {\rm F}oliation } (or CERF) if there 
    exists  a connected compact manifold $N$, endowed 
    with an SRF $\HH$, and a foliated imbedding $(M,\F) \subset (N,\HH)$
    such that $M$ is the regular stratum of $\SH$, that is, $M = R_{\HH}$. 
    We shall say that $(N,\HH)$ is a {\em zipper} of $(M,\F)$.  Both
    manifolds, $M$ and $N$, are connected or not at the same time.
    
    When $M$ is compact and $\F$ is regular then $(M,\F)$ is a CERF 
    and $(M,\F)$ the zipper. But in the general case, 
the zipper may
    not be unique.  The foliated manifold $(M, \F) = (\S^{3} \times 
    ]0,1[, \Ho\times \mathcal{I}))$, with $\Ho$ the Hopf foliation 
    of $\S^{3}$
    and $\I$ the foliation of $]0,1[$ by points,  is a CERF 
    possessing two natural zippers $(N_{i},\HH_{i})$, 
    $i=1,2$:
    \begin{itemize}
	\item[-] $N_{1} = \S^{4}$ and  $\HH_{1} $ is given by the orbits of 
	the
	$\sbat$-action: $z \cdot (z_{1},z_{2},t) =   (z \cdot z_{1}, z \cdot 
	z_{2},t) $, where $\S^{4} = \{ (z_{1},z_{2},t) \in \C \times \C 
	\times \R /  | / |z_{1}|^{2} + |z_{2}|^{2} +t^{2} = 1 \}$.
	\item[-]  $N_{2} =\C\P^{2}$ and  $\HH_{2} $ is given by the orbits of 
	the
	$\sbat$-action: $z \cdot [z_{1},z_{2},z_{3}] =   [z \cdot z_{1}, z \cdot 
	z_{2},z_{3}] $.
	\end{itemize}
	
    \medskip
    
    We consider in the sequel a manifold $M$ endowed with a CERF $\F$ and we fix $(N,\HH)$     a zipper.
    We present the Molino's desingularisation of $(N,\HH)$ in several 
    steps.

    \prg {\bf Foliated tubular neighborhood.} 
    A singular stratum $S \in \SH^{^{sin}}$ is a proper submanifold of the 
    riemannian manifold $(N,\HH,\mu)$.
    So it possesses a tubular neighborhood $(T_S,\tau_S,S)$.
    Recall that associated with this neighborhood we have the following
    smooth maps:
    \begin{itemize}
	\item[+] The {\em radius map} $\rho_S \colon T_S \to [0,1[$, which is
	defined fiberwise by $z\mapsto |z|$.  Each $t\not= 0$ is a regular
	value of the $\rho_S$.  The pre-image $\rho_S^{-1}(0)$ is $S$. 
	\item[+] The {\em contraction} $H_S \colon T_S \times [0,1] \to T_S$, which is 
	defined fiberwise by  $(z,r ) \mapsto r\cdot z$.  The restriction
	$(H_S)_t \colon T_S \to T_S$ is an imbedding for each $t\not= 0$ and
	$(H_{S})_0 \equiv \tau_S$.
    \end{itemize}

    \nt These maps verify $\rho_S(r \cdot u) = r \cdot \rho_S(u)$ 
    and $\tau_{S}(r \cdot u) = \tau_{S}(u)$. 
This tubular neighborhood can be chosen so that the two following 
important properties are verified (cf.  \cite{Mo}):

    \Zati Each  $(\rho_S^{-1}(t),\mathcal{H})$ is a 
	SRF, and

	 \zati Each  $(H_S)_{t} \colon  (T_S,\mathcal{F}) \to (T_S,\mathcal{F}) $ is a
	 foliated map.

	 \medskip

	 \nt 
When this happens, we shall say that $(T_S,\tau_S,S)$ is a {\em foliated tubular
neighborhood} of $S$.  The hypersurface $D_S = \rho_S^{-1}(1/2)$ is
the {\em core} of the tubular neighborhood.  We have $\depth \SHD <
\depth \SHTS$.

\medskip

There is a particular type of singular stratum we shall use in this work. A
stratum $S$ is a {\em boundary stratum} if $\codim_{N} \HH = 
\codim_{S}\HH_{S}- 1$.  The reason for this name is well illustrated
by the following example.  The usual $\sbat$-action on $\S^2$ by
rotations defines a singular riemannian foliation $\HH$ with two
singular leaves, two fixed points of the action.
These points are the boundary strata 
       and we have $
       N/\HH = [0,1]$. The boundary $\partial ( N/\HH )$ is given by 
       the boundary strata.
       In fact, the link of a boundary stratum is a sphere  with the 
       one leaf foliation (see, for example, \cite{SW} for the notion 
       of link).
       
\medskip
	 In the sequel, we shall use  the {\em partial blow up} 
$$\mathfrak{L}_S \colon (D_S\times [0,1[, \mathcal{H}\times 
\mathcal{I}) \to (T_S,\mathcal{H}),$$
which is the foliated smooth map defined by
$\mathfrak{L}_S(z,t) =2t \cdot z$. Here, $\mathcal{I}$ denotes the 
pointwise foliation.
The restriction  
\begin{equation}
    \label{form}
    \mathfrak{L}_{S} \colon (D_S \times ]0,1[ , \mathcal{H} \times 
\mathcal{I}) \to (T_S\menos S,\mathcal{H})
\end{equation}
is a foliated diffeomorphism.

\prg {\bf Foliated Thom-Mather system}.  In the proof of Lemma
\ref{reduc} we find two strata $S_{1} \preceq S_{2}$ endowed with two
tubular neighborhoods $T_{S_{1}}$ and $T_{S_{2}}$.  
We shall need $T_{S_{2}}\menos T_{S_{1}}$ to be a tubular neighborhood  
of $S_{2}\menos T_{S_{1}}$, but this is not always achieved. To guarantee this property, 
we  introduce the following notion, which is inspired in the abstract stratified
objects of \cite{Ma,Th}.

A family of foliated tubular neighborhoods $\{T_{S } \ | \ S \in
\SF^{^{sin}}\}$ is a {\em foliated Thom-Mather system} of $(N,\HH)$
if, for each pair of singular strata $S_{1}, S_{2}$ with $S_{1}
\preceq S_{2}$, we have
\begin{equation}
   \label{TM1}
   \rho_{S_{1}} = \rho_{S_{1}} \rondp \tau_{S_{2}} \ \ \hbox{ on
   } \ \ T_{S_{1}} \cap T_{S_{2}} = \tau_{S_{2}}^{-1}(T_{S_{1}} \cap
   S_{2}).
\end{equation} 
    In these conditions we have the property:
    \begin{equation}
       \label{TM2}
	\rho_{S_{1}}\underbrace{(r \cdot z)}_{H_{S_{2}}(z,r)}= \rho_{S_{1}}(z),
	\end{equation}
	for each $r \in [0,1]$ and $u \in T_{S_{1}} \cap T_{S_{2}}$.
	We conclude that the restriction
\begin{equation}
    \label{menos}
	\tau_{S_{2}} \colon \bigg( T_{S_{2}}\menos \rho_{S_{1}} ^{-1}(I) 
	\equiv \tau_{S_{2}}^{-1}(S_{2}\menos \rho_{S_{1}} ^{-1}(I) )
	\bigg)
	\TO 
	\left( S_{2}\menos \rho_{S_{1}}^{-1}(I) \right)
\end{equation}
	is a foliated tubular neighborhood of $S_{2}\menos \rho^{-1}(I)$ on 
	$N\menos \rho_{S_{1}}^{-1}(I)$, where $I \subset [0,1[$ is a closed subset.
	The foliated diffeomorphism \refp{form} becomes
	\begin{equation}
	    \label{form2}
	    \mathfrak{L}_{S_{2}} 
	    \colon (D_{S_{2}}\menos \rho_{S_{1}}^{-1}(I))\times ]0,1[ , 
	    \mathcal{H} \times 
	\mathcal{I}) \to ((T_{S_{2}}\menos {S_{2}})\menos 
	\rho_{S_{1}}^{-1}(I),
	\mathcal{H})
	\end{equation}
    \bP
    \label{TM}
    Each compact manifold endowed with an SRF possesses a foliated
    Thom-Mather system.  \eP \pro See Appendix.  \qed \prg {\bf Blow up}. 
    The Molino's blow up of the foliation $\mathcal{H}$ produces a new
    foliation $\widehat{\mathcal{H}}$ of the same kind but with smaller
    depth (see \cite{Mo0} and also \cite{SW2}).  The main idea is to
    replace each point of the closed strata by a sphere.

\smallskip

We suppose that $\depth \SH > 0$.  The union of closed (minimal) strata
we denote by $ \Sigma_0$.  We choose $T_0$ a disjoint family of
foliated tubular neighborhoods of the closed strata.  The union of the
associated cores is denoted by $D_0$.  Let $\mathfrak{L}_0 \colon (D_0
\times [0,1[, \HH \times \mathcal{I}) \to (T_0 ,\HH)$ be the
associated partial blow up.  The {\em blow up} of $(N,\HH,\mu)$ is
 $$
 \mathfrak{L} \colon ( \wh{N},\wh{\HH},\wh{\mu})  \TO (N,\HH,\mu)
 $$
 where
 \begin{itemize}
     \item[-] The manifold $\wh{N}$ is 
$$
\wh{N} = 
\left\{
\Big( D_0 \times ]-1,1[\Big) \coprod \Big( (N\menos \Sigma_{0}) \times \{ 
-1,1\}\Big)
\right\} \Big/ \sim ,
$$
where  $(z,t) \sim (\mathfrak{L}_0 (z, |t|) , t/|t|)$.
 Notice that $D_0 \times ]-1,1[$ and $(N\menos \Sigma_{0}) \times \{ 
 -1,1\}$ are open subsets of $\wh{N}$ with 
 $$
 (D_0 \times ]-1,1[) \cap ((N\menos \Sigma_{0}) \times \{
 -1,1\}) = 
 D_0 \times (]-1,0[ \cup ]0,1[).
 $$
%
 \item[-] The foliation $\wh{\HH}$ is determined by 
 $$
 \wh{\HH}|_{D_0 \times ]-1,1[} = 
 \HH|_{{D_0} }\times \mathcal{I} \ \hbox{ and } \ 
 \wh{\HH}|_{(N \menos  \Sigma_{0}) \times \{ -1,1\} } = 
 \HH|_{{N \menos  \Sigma_{0} } }\times \mathcal{I} .
 $$
 Here, $\mathcal{I}$ denotes the 0-dimensional foliation of $]-1,1[$.
 
 \item[-] The riemannian metric $\wh\mu$ is
 $$
f\cdot (\mu|_{D_0} + dt^{2}) + (1- f ) \cdot \mu|_{N\menos \Sigma_{0}},
 $$
 where $ f\colon \wh{N} \to [0,1]$ is the smooth map defined by
 $$
f(v) = \left\{
 \begin{array}{ll}
     \xi (|t|) & \hbox{if } v =(z,t) \in D_0 \times ]-1,1[\\[,2cm]
      0 & \hbox{if } v =(z,j) \in(N  \menos   \rho^{-1}_0([0,3/4]) \times \{-1,1\},
     \end{array}
     \right.
 $$
 with $\xi \colon [0,1] \to [0,1]$ a smooth map verifying $\xi
 \equiv 1 $ on $[0,1/4]$ and $\xi\equiv 0$ on $[3/4,1[$.
 
 \item[-] The map $\mathfrak{L}$ is defined by
$$
\mathfrak{L}(v) = 
\left\{
\begin{array}{ll}
    \mathfrak{L}_0 (z, |t|)  & \hbox{if } v = (z,t) \in  D_0 \times 
]-1,1[ \\[,2cm]
z& \hbox{if }  v = (z,j) \in (N\menos  \Sigma_{0}) \times \{ 
-1,1\},
\end{array}
\right.
$$
\end{itemize}
Notice that the blow up of $(N,\HH,\mu)$ depends just on the choice 
of $\xi$.  So, we fix from now on  such a $\xi$.

\prgg {\bf Remarks.}

\Zati The blow up of 
$(N \times \R,\HH \times \I ,\mu + dt^{2})$ is just 
$
\mathfrak{L} \times \Ide \colon  (\wh{N} \times \R,\wh{\HH}\times 
\I,\wh{\mu} + dt^{2})  \TO (N \times \R,\HH \times \I ,\mu + dt^{2}).
$

\zati The manifold $\wh{N}$ is connected and compact, the foliation  $\wh{\HH}$ is an SRF
and $\wh{\mu}$ is an adapted metric. 

\zati The map $\mathfrak{L}$ is a foliated continuous map whose restriction 
to 
$(\wh{N} \menos  \mathfrak{L}^{-1}(\Sigma_{0}),\wh{\HH}) 
\equiv ((N\menos \Sigma_{0} )\times \{ -1,1\}, \HH \times 
\mathcal{I})$ is the canonical projection on the first factor.

\zati We shall denote by $M_{1}$ the regular stratum $R_{\wh{\HH}}$ and $\F_{1}$ the restriction of 
$\wh{\HH}$  to $M_{1}$. In fact, the foliation $\F_{1}$ 
 is a CERF on $M_{1}$ and a zipper is given by $(\wh{N}, \wh{\HH})$.
Notice that we have the inclusion 
$
\mathfrak{L}^{-1} (M) = M \times \{ -1,1\}\subset  M_{1}.
$
We choose  $\mathcal{S}_{1} \colon M\to M_{1}$ a 
smooth foliated imbedding with $\mathfrak{L} \rondp \mathcal{S}_{1} = 
\Ide$. 
There are two of them.

\zati {\em The stratification $\SHhat$}. For each non-minimal stratum $S \in \SH$ there exists a unique 
stratum 
$S^{\ \widehat{}} \in \SHhat$, with $\mathfrak{L}^{-1}(S) \subset S^{\ \widehat{}}  $,
in fact,
$
S^{\ \widehat{}} = 
\left\{
\Big(\left( D_0  \cap S\right) \times ]-1,1[\Big) \coprod \Big( S\times \{ 
-1,1\}\Big)
\right\} \Big/ \sim .
$
This gives 
$
    \SHhat = \{ S^{\ \widehat{}}  \ | \ S \in \SH \hbox{ and  
    non-minimal}\}.
$
We have the following important properties
\begin{itemize}
    \item[-] $\depth_{\HH} S -1 = \depth_{\wh{\HH}}S^{\ 
    \widehat{}} $, for each non-minimal stratum $S \in \SH$,
    \item[-] $\mathfrak{L}^{-1}(\Sigma_{i}\menos \Sigma_{0}) 
    =\Sigma_{i-1}(\wh{N})\menos 
    \mathfrak{L}^{-1}(\Sigma_{0}) )$ for each $i 
    \in \Z$,  and
    \item[-] 
$
\depth \SHhat < \depth \SH.
$
\end{itemize}

\zati We shall use the diffeomorphism $\sigma \colon \wh{N} \to \wh{N}$ 
defined by
$$
\sigma(v) = 
\left\{
\begin{array}{ll}
    (z,-t) & \hbox{if } v = (z,t) \in  D_0 \times 
]-1,1[ \\[,2cm]
(z,-j)& \hbox{if }  v = (z,j) \in (N\menos  \Sigma_{0}) \times \{ 
-1,1\}.
\end{array}
\right.
$$
In fact, the diffeomorphism $\sigma$ is a foliated isometry verifying 
$\mathfrak{L} \rondp \sigma = \mathfrak{L}$. 
It induces the smooth foliated action $\Phi \colon \Z_{2}\times M_{1} \to M_{1}$ 
defined by
$
\zeta \cdot v = \sigma(v),
$
where $\zeta$ is the generator of $\Z_{2}$.

%

\prg {\bf Molino's desingularisation}. If the depth of $\SHhat$ is not 
0 then the blow up can be continued (cf. 1.6.1 (b)). 
In the end, we have a riemannian foliated manifold 
$(\widetilde{N},\widetilde{\mathcal{H}},\wt\mu)$ and a foliated continuous map 
$\mathfrak{N} \colon (\widetilde{N} ,\wt\HH)\to (N,\HH)$, whose restriction 
$
\mathfrak{N} \colon \mathfrak{N}^{-1}(M)
\to M
$
is a smooth trivial bundle (cf. 1.6.1 (d)).
Notice that $\wt{N}$ is connected and compact.
This type of construction is a {\em Molino's desingularisation} of 
$(N,\HH,\mu)$ (cf.  \cite{Mo0}).

We choose $\mathcal{S}\colon M\to \wt{N}$ a 
smooth foliated imbedding verifying $\mathfrak{N} \rondp \mathcal{S} 
=\Ide$. It always exists.

    \section{Tautness}
    The tautness of  a RF on a compact manifold can be detected 
    by using the basic cohomology. This is not the case 
    when the manifold is not compact. In this section we recover this result 
    for  a CERF.
    
    \medskip
    
    We consider in the sequel a manifold $M$ endowed with a CERF 
    $\F$ and we fix  a zipper 
    $(N,\HH)$. We also consider a 
  Molino's desingularisation $\mathfrak{N} \colon (\wt{N},\wt\HH) \to 
  (N,\HH)$. We shall write $n = \codim \F$. We fix  $\mathcal{S}\colon M\to \wt{N}$ a 
smooth foliated imbedding  verifying $\mathfrak{N} \rondp \mathcal{S} 
=\Ide$.

    \prg {\bf Basic cohomology}.  Recall that the {\em basic cohomology}
    $\hiru{H}{*}{\mf}$ is the cohomology of the complex
    $\hiru{\Om}{*}{\mf}$ of {\em basic forms}.  A differential form $\om$
    is basic when $i_{X}\om = i_{X}d\om =0$ for every vector field $X$
    tangent to $\F$.
    
    An open covering $\{ U , V \}$  of $M$ by saturated open 
    subsets possesses a 
    subordinated partition of the unity made up of  basic  
    functions (see Lemma below). 
    For such a covering we have the Mayer-Vietoris short sequence
  \begin{equation}
      \label{mv}
    0 \to \hiru{\Om}{*}{\mf} \to 
    \hiru{\Om}{*}{U/\mathcal{F}} \oplus 
    \hiru{\Om}{*}{V/\mathcal{F}}  \to
    \hiru{\Om}{*}{(U \cap V)/\mathcal{F}}  \to 0,
    \end{equation}
    where the maps are defined by restriction. The third map is onto since the 
    elements of the partition of the unity are basic  
    functions. Thus,  
    the sequence
     is exact.

     The {\em compactly supported  basic cohomology} $\lau{H}{*}{c}{\mf}$
     is the 
     cohomology of the basic subcomplex  $\lau{\Om}{*}{c}{\mf}
     = \{ \om \in \hiru{\Om}{*}{\mf} \ |\  \hbox{ the support of 
     $\om$ is compact}\}$.

     The {\em twisted basic cohomology} $\lau{H}{*}{\kappa}{\mf}$, 
     relatively to the cycle $\kappa \in \hiru{\Om}{1}{\mf}$, is the 
     cohomology of the basic complex $\hiru{\Om}{*}{\mf}$ relatively 
     to the differential $ \om \mapsto d\om - \kappa \wedge \om$.  This
     cohomology does not depend on the choice of the cycle: we have
     $\lau{H}{*}{\kappa}{\mf}\cong \lau{H}{*}{\kappa + df}{\mf}$ through
     the isomorphism: $[\omega] \mapsto [e^f \omega]$.
     
     Given $V$, a $\Z_{2}$-invariant saturated open subset of $M$, we 
     shall write
     \begin{eqnarray*}
\left( \hiru{H}{*}{V/\F} \right)^{\Z_{2}} &= &\{ \omb \in 
\hiru{H}{*}{V/\F}  \ | \ \sigma^{*}\omb = \omb \}\\
\left( \hiru{H}{*}{V/\F} \right)^{-\Z_{2}} &=& \{ \omb \in 
\hiru{H}{*}{V/\F}  \ | \ \sigma^{*}\omb = -\omb \}.
	 \end{eqnarray*}
\medskip

For the existence of the Mayer-Vietoris sequence \refp{mv} we need the 
following folk result, well-known for compact Lie group actions and 
regular riemannian foliations.
\bL
\label{sat}
Any  covering   of $M$ by saturated open 
subsets possesses a 
subordinated partition of the unity made up of  basic  
functions.
\eL
\pro
The closure $\overline{L}$ of a leaf $L \in \F$ is a saturated 
submanifold of $M$ whose leaves are dense (cf. \cite{Mo}). So, the open subsets $U$ 
and $V$ are in fact $\overline{\F}$-saturated subsets.

The  closure $\overline{L}$ of a leaf $L \in \F$ possesses a 
 tubular neighborhood as in 1.4 (cf. \cite{Mo}). 
Since the family of these tubular neighborhoods is a basis for the 
family of saturated open subsets then it suffices to construct the 
partition of unity relatively to the tubular neighborhoods. This is 
done by using the radius maps.
\qed

 \prg {\bf Tautness reminder (compact case)}.  Given a bundle-like metric
 $\mu$ on $(M,\F)$, the mean curvature form $\kappa_{\mu} \in
 \hiru{\Om}{1}{M}$ is defined as follows (see for example \cite{To}). 
 Consider the second fundamental form of the leaves and $W$ the
 corresponding Weingarten map.  Then, $$\kappa_{\mu}(X) = \left\{
    \begin{array}{ll}
	\trace W(X) & \hbox{if $X$ is orthogonal to the foliation $\F$}\\[,3cm]
0& \hbox{ if
    $X$ is tangent to the foliation $\F$.}
    \end{array}
    \right.
    $$
    When the manifold is compact (and then $(N,\HH) = (M,\F)$), the following 
    properties of $\kappa_{\mu}$  are well-known:
    
    \Zati  The form $\kappa_{\mu}$ 
    can be supposed to be basic, i.e., there exists a a bundle-like 
    metric
    $\mu$ such that its mean curvature form is basic (see \cite{D}). 
    
    \zati  If $\kappa_{\mu}$ is basic, then $\kappa_{\mu}$ is a closed form 
    (see \cite{To}). The Example 2.4 of
   \cite{CE} shows that the compactness assumption cannot be removed: 
   there, the mean curvature form is basic, but not closed.
   
    \zati The class $\kappab = [\kappa_{\mu}] \in \hiru{H}{1}{\mf}$ 
    does not depend on the metric, 
    but just on $\F$ (see \cite{Suso}). This is the {\em tautness 
    class} of $\F$.
    
    The mean curvature form contains some  geometric information about $\F$. 
    Recall that the
    foliation  $\F$ is {\em taut} if there exists a riemannian metric $\mu$ 
    on $M$ such that every leaf is a
    minimal submanifold of $M$. It is known (see \cite{To}) that 
    $$
  \F \hbox{ is 
    taut } \Longleftrightarrow  \hbox{the tautness class
    $\kappab$ vanishes}. 
    $$
    We also have the following cohomological characterizations for the 
    tautness of $\F$:
    $$
  \F \hbox{ is 
    taut } \Longleftrightarrow  \hiru{H}{n}{\mf} \ne 0,
    $$
    when $\F$ is transversally oriented.  We also have that
    $$
   \F \hbox{ is 
     taut } \Longleftrightarrow  \lau{H}{0}{\kappa_{\mu}}{\mf} \ne 0,
     $$
     when $M$ is oriented and $\F$ is transversally oriented.
    
     Immediate examples of taut foliations  are
    isometric flows (i.e. 1-foliations induced by the orbits of a nonvanishing 
    Killing vector
    field), isometric actions on compact manifolds  and compact foliations with locally 
    bounded volume of leaves (foliations where every leaf is compact, 
    see \cite{EMS},\cite{RU}). 
    
    In the example of \cite{CE} referred above, we find a non-compact
    manifold where the tautness class may not exist.
    
    These results are not directly extendable to the framework of singular
    riemannian foliations, as we can see in the following example: the usual
    $\sbat$-action on $\S^2$ by rotations defines a singular riemannian
    foliation $\HH$ with two singular leaves, two fixed points.  Notice that
    $\hiru{H}{1}{\nh} = 0$.  But there cannot exist a metric on $\S^2$ such
    that the one dimensional orbits are geodesics (in dimension 1,
    ÒminimalÓ implies ÒgeodesicÓ).  To see this, it suffices to consider a
    totally convex neighborhood $U$ of one of the fixed points: it would
    contain as geodesics some orbits (full circles, indeed), apart from
    rays, which would contradict uniqueness of geodesics connecting any
    two points of $U$.  In higher dimensions (i.e. the cone of a sphere),
    one may consider the volume of every leaf.  It should be constant by
    minimality, which would give a positive volume to the vertex (for more
    details see \cite{MW}).

\prg {\bf Construction of $\kappab$}.
We would like to define a cohomological class $\kappab \in 
\hiru{H}{1}{M/\F}$ which would play a r™le similar to that of  the tautness 
class of a regular riemannian foliation on a regular manifold. First,
we introduce the notion of D-metric.
A bundle-like metric metric $\mu$ on $(M,\F)$ is a {\em D-metric} if the mean 
curvature form $\kappa_{\mu}$ is a basic cycle.  The {\em tautness 
class} of $\F$  is the cohomological class $\kappab
= [\kappa_\mu] \in \hiru{H}{1}{M/\F}$. 
The next Proposition proves that this class is well defined and 
independent of the D-metric.

\bL
\label{reduc}
There exists a saturated open subset $U  \subset M$  such that

\Zati the inclusion $\upsilon \colon U \hookrightarrow M$ induces the 
isomorphism $\upsilon^{*} \colon \hiru{H}{*}{M/\F} \to 
\hiru{H}{*}{U/{\F}}$,  and

\zati the closure $\overline{U}$ (in $N$) is included in $M$.

\eL
\pro 
We consider $\{T_{S} \ | \ S \in \SF\}$ a foliated Thom-Mather system 
of $(N,\HH)$ (cf. Proposition \ref{TM}). For each $i\in \Z$ we 
write: 
\begin{itemize}
\item[-] $\tau_{i}\colon T_{i} \to \Sigma_{i}\menos \Sigma_{i-1}$ 
the associated foliated  tubular neighborhood of $\Sigma_{i}\menos \Sigma_{i-1} $.
\item[-] $\rho_{i} \colon T_{i}\to [0,1[$ its radius function, and
\item[-] $D_{i}$ the core of $T_{i}$.
\end{itemize}

The family 
$\{M \cap T_{0},M\menos \rho_{0}^{-1}([0,7/8])\}$ is a saturated 
open covering of $M$.
Notice that the inclusion $((M \cap T_{0} )\menos \rho_{0}^{-1}([0,7/8]),\F) 
\hookrightarrow (M\cap T_{0},\F) $ induces an 
isomorphism for the basic cohomology since it is foliated 
diffeomorphic to the inclusion
$$
((M \cap D_{0})\times ]7/8,1[ , \F \times \I) \hookrightarrow 
((M \cap D_{0} )\times ]0,1[ , \F \times \I)
$$
(cf. \refp{form}).
From the Mayer-Vietoris sequence we conclude that the inclusion 
$
M\menos \rho_{0}^{-1}([0,7/8])
\hookrightarrow 
M
$
induces the isomorphism
$$
\hiru{H}{*}{M/\F} =
\hiru{H}{*}{M\menos \rho_{0}^{-1}([0,7/8])/\F}
$$
(cf. \refp{mv}).

Notice now that 
$M\menos \rho_{0}^{-1}([0,7/8])$ is the regular part 
of $(N\menos\rho_{0}^{-1}([0,7/8]),\HH)$. Moreover, the family
$$
\{ T_{S} \menos  \rho_{0}^{-1}([0,7/8]) \ | \ S \in \SF  ,   \dim S > 0\}
$$
is a foliated Thom-Mather system of
$(N\menos\rho_{0}^{-1}([0,7/8]),\HH)$ (cf.  \refp{menos}).  The same
previous argument (using the foliated diffeomorphism \refp{form2}
instead of that of \refp{form}) gives
$$
\hiru{H}{*}{M\menos \rho_{0}^{-1}([0,7/8]) /\F}
=
\hiru{H}{*}{M\menos (\rho_{0}^{-1}([0,7/8]) \cup \rho_{1}^{-1}([0,7/8]) )/\F}.
$$
So, one gets the isomorphisms
$$
\hiru{H}{*}{M/\F} 
=
\cdots 
=
\hiru{H}{*}{M\menos (\rho_{0}^{-1}([0,7/8]) \cup  \cdots \cup 
\rho_{p-1}^{-1}([0,7/8]) )/\F},
$$
where $p = \depth \SH$.
Take $U = M\menos (\rho_{0}^{-1}([0,7/8]) \cup  \cdots \cup 
\rho_{p-1}^{-1}([0,7/8]) )$, which is an open saturated subset of 
$\F$ included on $M$. This gives (a).

Consider $K = M\menos (\rho_{0}^{-1}([0,6/8[) \cup  \cdots \cup 
\rho_{p-1}^{-1}([0,6/8[) )$, which is a subset of $M$ containing $U$.
We compute its closure on $N$:
\begin{eqnarray*}
    \overline{K} &=& \overline{M\menos (\rho_{0}^{-1}([0,6/8]) \cup  \cdots \cup 
\rho_{p-1}^{-1}([0,6/8]) )} \subset 
\overline{M}\menos \left((\rho_{0}^{-1}([0,6/8[))^{\hbox{¡}} \cup  \cdots \cup 
(\rho_{p-1}^{-1}([0,6/8[))^{\hbox{¡}}\right)\\
&=& N\menos \left(\rho_{0}^{-1}([0,6/8[) \cup  \cdots \cup 
\rho_{p-1}^{-1}([0,6/8[)\right) = M\menos (\rho_{0}^{-1}([0,6/8[) \cup  \cdots \cup 
\rho_{p-1}^{-1}([0,6/8[) ),
\end{eqnarray*}
since $N\menos M = \Sigma_{p-1} = 
\rho_{0}^{-1}(\{0 \}) \cup  \cdots \cup 
\rho_{p-1}^{-1}(\{ 0 \}) $.
This implies that $K$ is a closed subset of $N$ and therefore 
compact. This gives (b).
\qed
 \bP
 \label{existe}
The tautness class of 
a CERF exists and it does not depend on the choice of the 
D-metric.
\eP
\pro 
We proceed in two steps.

\smallskip

{\em (i) - \underline{Existence of D-metrics}.} Since $\wt{N}$ is compact then 
there exists a  D-metric $\nu$ on $\wt{N}$. 
The metric  $\mathcal{S}^{*}\nu$ is a bundle-like metric on $M$. 
Since $\kappa_{\mathcal{S}^{*}\nu} = \mathcal{S}^{*}\kappa_{\nu}$ then 
$\mathcal{S}^{*}\nu$ is a D-metric on $M$. 

\smallskip

{\em (ii) - \underline{Uniqueness  of $\kappab$}}.
Consider $\mu$ a D-metric on $M$ and let  $\kappabt$  be the tautness class
of $\wt\HH$.  It suffices to prove $ [\kappa_{\mu}] =
\mathcal{S}^{*}\kappabt.  $ Take $U$ as in the previous Lemma and
$\{f,g\}$ a subordinated partition of unity associated to the covering
$\{M, N\menos \overline{U}\}$ made up of basic functions (cf.  Lemma
\ref{sat}).  Notice that $f \equiv  1$ on $U$.  So, it suffices to prove
\begin{equation}
    \label{kapa}
    \upsilon^{*}[\kappa_{\mu}] = \upsilon^{*}\mathcal{S}^{*}\kappabt 
    \end{equation}
    (cf Lemma  \ref{reduc} (a)). Consider the following diagram:
    $$ 
    \begin{picture}(200,70)(0,-10) 
    \put(0,0){\makebox(0,0){$U$}}
    \put(100,0){\makebox(0,0){$M$}} 
    \put(200,0){\makebox(0,0){$N$,}} 
    \put(0,50){\makebox(0,0){$\mathfrak{N}^{-1}(U)$}}
    \put(100,50){\makebox(0,0){$\mathfrak{N}^{-1}(M)$}}
    \put(200,50){\makebox(0,0){$\wt{N}$}}

    \put(20,0){\vector(1,0){54}}
    \put(120,0){\vector(1,0){54}}
    \put(25,50){\vector(1,0){50}}
    \put(125,50){\vector(1,0){50}}
    \put(125,8){\vector(2,1){58}}
    
    \put(-5,40){\vector(0,-1){30}}
    \put(5,10){\vector(0,1){30}}
   \put(95,40){\vector(0,-1){30}}
       \put(105,10){\vector(0,1){30}}
       \put(200,40){\vector(0,-1){30}}

       \put(47,8){\makebox(0,0){$\upsilon$}}
       \put(50,58){\makebox(0,0){$\upsilon'$}}
       \put(150,58){\makebox(0,0){$\iota$}}
       \put(150,8){\makebox(0,0){$\iota'$}}
       
       \put(-15,25){\makebox(0,0){$\mathfrak{N}''$}}
    \put(15,25){\makebox(0,0){$\mathcal{S}''$}}
    \put(85,25){\makebox(0,0){$\mathfrak{N}'$}}
    \put(115,25){\makebox(0,0){$\mathcal{S}'$}}
    \put(145,25){\makebox(0,0){$\mathcal{S}$}}
    \put(190,25){\makebox(0,0){$\mathfrak{N}$}}
    \end{picture}                                            
    $$
   where $\upsilon'$, $\iota$, $\iota'$ are the natural inclusions, 
   $\mathfrak{N}'$, $\mathfrak{N}''$ are the 
   restrictions of $\mathfrak{N}$ and $\mathcal{S}'$, $\mathcal{S}''$ 
   are the restrictions
   of $\mathcal{S}$. Notice the equalities
	$
\mathcal{S} \rondp \upsilon = 
	\iota \rondp \upsilon'
	\rondp \mathcal{S}''
	$
and
$
\upsilon = \mathfrak{N}'\rondp \upsilon'
	\rondp \mathcal{S}''.
	$

	The maps $\upsilon'$, ${\mathfrak{N}'}$ are foliated local diffeomorphism so $
	\kappa_{\upsilon'^*{\mathfrak{N}'}^*\mu} =
	\upsilon'^*\mathfrak{N}'^*\kappa_\mu.  $ This differential form is a
	cycle since $\mu$ is a D-metric.  Take on $\wt{N}$ the riemannian
	metric $${\lambda} = \mathfrak{N}^{*}f \cdot \mathfrak{N}'^{*} \mu +
(1-\mathfrak{N}^{*}f ) \cdot \nu$$
It is a bundle-like metric since $f$ 
is basic and the support of $f$  is included on $M$. 
Notice the equality:
$
\upsilon'^{*}\iota^{*}\lambda = \upsilon'^{*}\mathfrak{N}'^{*}\mu.
$

We can use this metric for the computation of  $\kappabt$ in the 
following way.
Let $(\kappa_{\lambda})_{_{b}}$ be  the basic part of the mean curvature form 
$\kappa_{\lambda}$,
  relatively to the $\lambda$-orthogonal decomposition
 $\hiru{\Om}{*}{\wt{N}}=\hiru{\Om}{*}{\widetilde{N}/\widetilde{\HH}} \oplus
 \hiru{\Om}{*}{\widetilde{N}/\widetilde{\HH}}^{\bot}$. It is a basic cycle and we have 
 $\kappabt = [(\kappa_{\lambda})_{_{b}}]$ (see \cite{Suso}).

Since $\mathfrak{N}^{-1}(U)$ is $\wt{\HH}$-saturated open subset of $
\wt{N}$ then it is also  a $\overline{\wt{\HH}}$-saturated subset  
 of $\wt{N}$ (see the proof of Lemma \ref{sat}). 
 From the definition of the basic component $(\kappa_{\lambda})_{_{b}}$ we
 get that the restriction
 $\upsilon'^{*}\iota^{*}(\kappa_{\lambda})_{_{b}}$ is defined by using just
 $(\mathfrak{N}^{-1}(U),\upsilon'^{*}\iota^{*}\wt\HH ,
 \upsilon'^{*}\iota^{*}\lambda)$.  As a consequence, we have
 $\upsilon'^{*}\iota^{*}\left((\kappa_{\lambda})_{_{b}}\right)=
 \kappa_{\upsilon'^{*}\iota^{*}\lambda}= \kappa_{
 \upsilon'^*\mathfrak{N}'^*\mu} =
 \upsilon'^*\mathfrak{N}'^*\kappa_\mu$.
 
 Finally, we get
$$
 \upsilon^{*}\mathcal{S}^{*}\kappabt 
=
\mathcal{S}''^{*} \upsilon'^{*}\iota^{*}\kappabt
=
\mathcal{S}''^{*} \upsilon'^{*}\iota^{*}[(\kappa_{\lambda})_{_{b}}]
=
 \mathcal{S}''^{*}[\upsilon'^*\mathfrak{N}'^*\kappa_\mu]
=
 \mathcal{S}''^{*}\upsilon'^{*}\mathfrak{N}'^{*}[\kappa_\mu]
 =
 \upsilon^{*}[\kappa_\mu].
 $$
This gives \refp{kapa}.
\qed

\prgg {\bf Remarks.} 

\Zati The tautness class $\kappab$ of $(M,\F)$ and the tautness class 
$\kappabt$ of a Molino's 
desingularisation $(\wt{N},\wt\HH)$ are 
related 
by the formula
$
\kappab = \mathcal{S}^{*}\kappabt,$
where $\mathcal{S}   \colon M\to 
\wt{N}$ is any smooth foliated
embedding verifying $\mathfrak{N}\rondp \mathcal{S} = \Ide$.

%
 \zati 
 Let  $\kappab_{1}$  be the tautness class of $(M_{1},\F_{1})$.  This class is
 $\Z_{2}$-invariant (cf.  1.6).  This comes from the fact that the
 diffeomorphism $\sigma$ preserves $\HH_{1}$ and therefore
 $\kappa_{\sigma^*\mu} = \sigma^{*}\kappa_{\mu}$ for a D-metric $\mu$
 on $M_{1}$.  Thus, the metric $\sigma^{*}\mu$ is also a D-metric and
 we obtain: $ \sigma^*\kappab_{1} = [\sigma^{*}\kappa_{\mu}] =
 [\kappa_{\sigma^{*}{\mu}}] = \kappab_{1}.  $

\prg{\bf First characterization of tautness: vanishing of $\kappab$.}

We give the first characterization of the tautness of $\F$ through 
the vanishing of $\kappab$. We lift the question to the Molino's 
desingularisation $(\wt{N},\wt\HH)$.

\bL
The map $\mathcal{S}_{1}\colon M
\to M_{1}$ (cf. 1.6)
induces the isomorphism:
$$
\left( \hiru{H}{*}{M_{1}/\mathcal{F}_{1}}\right)^{\Z_{2}}
\cong
\hiru{H}{*}{M/\mathcal{F}}.
$$
\eL
\pro 
The open covering $\left\{ (D_0\cap M) \times ]-1,1[, 
M \times \{ -1,1\}\right\}$ of $M_{1}$ is  a $\Z_{2}$-equivariant 
one (cf. 1.6.1 (f)). 
So, from Mayer-Vietoris (cf. 
\refp{mv}) we get 
the long exact sequence
\begin{eqnarray*}
    \! \! \! \!  \cdots \! \! &\to&  
   \left( \hiru{H}{j}{(D_0\cap M) \times ]-1,1[/\F \times 
   \mathcal{I}}\right)^{\Z_{2}}
   \oplus 
   \left( \hiru{H}{j}{M \times \{ -1,1\}/\F \times 
   \mathcal{I}} \right)^{\Z_2}
  \stackrel{I}{\TO}
   \\
   &\to &    \left( \hiru{H}{j}{(D_0\cap M) \times 
   (]-1,0[ \cup ]0,1[)/\F \times 
   \mathcal{I}}\right)^{\Z_2}
   \to
   \left( \hiru{H}{j+1}{M_{1}/\F_{1}}\right)^{\Z_2}
   \stackrel{I}{\TO}
    \\
    &\to &  \left( \hiru{H}{j+1}{(D_0\cap M) \times ]-1,1[/\F \times 
      \mathcal{I}}\right)^{\Z_2}
      \oplus 
      \left( \hiru{H}{j+1}{M \times \{ -1,1\}/\F \times 
      \mathcal{I}} \right)^{\Z_2} \to \cdots ,
\end{eqnarray*}
where $I$ denotes the restriction map (i.e. induced by the inclusion).
Since the natural projection $(D_0\cap M) \times 
]-1,1[ \to (D_0\cap M) $ is 
$\Z_{2}$-invariant, then we get isomorphisms
\begin{eqnarray*}
    \hiru{H}{*}{(D_0\cap M) /\F 
       }&\cong &
       \left( \hiru{H}{*}{(D_0\cap M) \times ]-1,1[/\F \times 
   \mathcal{I}}\right)^{\Z_{2}}
\\
\hiru{H}{*}{(D_0\cap M) /\F 
      }  &\cong&
      \left( \hiru{H}{*}{(D_0\cap M) \times 
	 (]-1,0[ \cup ]0,1[)/\F \times 
	 \mathcal{I}}\right)^{\Z_{2}}.
	 \end{eqnarray*}
We conclude that the inclusion $M \times \{-1,1\}\hookrightarrow M_{1}$ induces 
the isomorphism
$$
   \left( \hiru{H}{*}{M_{1}/\F_{1}}\right)^{\Z_{2}} \cong  
\left( \hiru{H}{*}{M \times \{ -1,1\}/\F \times 
      \mathcal{I}}\right)^{\Z_{2}}.
   $$
   Since the natural projection $P \colon M \times \{ -1,1\} 
   \to
 M$ is 
   $\Z_{2}$-invariant, then we get isomorphism
   $$
\hiru{H}{*}{M /\F } 
	\cong
	\left( \hiru{H}{*}{M \times \{ -1,1\}/\F \times 
	     \mathcal{I}} \right)^{\Z_{2}}.
	     $$
	     Since $P \rondp\mathcal{S}_{1} = \Ide$ then  
	     $\mathcal{S}_{1}$ induces the 
	     isomorphism 
	     $
	      \left( \hiru{H}{*}{M _{1}/\F} \right)^{\Z_{2}}
		   \cong
		\hiru{H}{*}{M /\F} .
			$
		 \qed

The main result of this section is the following.
\bT
\label{T1}
Let $M$ be a manifold endowed with a CERF $\F$. 
Then, the following two  statements are equivalent:

\Zati The foliation $\F$ is taut.

\zati The tautness class $\kappab \in \hiru{H}{1}{M/\F}$ vanishes.
\eT
\pro
We prove the two implications.

\medskip

$(a) \Rightarrow (b)$. There exists a D-metric $\mu$ on $M$ 
with $\kappa_{\mu}= 0$. Then $\kappab = [\kappa_{\mu}] = 0$.

\medskip

$(b) \Rightarrow (a)$. We proceed by induction on $\depth \SH^{^{sin}}$. When 
this depth is 0 we have the regular case of \cite{Suso}.
When this depth is not 0 then we can consider the Molino's 
desingularisation
$
\mathfrak{N}= \mathfrak{L} \rond \mathfrak{N}_{1} \colon 
(\wt{N},\wt{\HH})  
\to (N,\HH)$
of $(N,\HH)$, where 
$
\mathfrak{N}_{1}  \colon (\wt{N},\wt{\HH}) 
\to (\wh{N},\wh{\HH})$
is a Molino's desingularisation of $(\wh{N},\wh{\HH})$ (cf. 1.6 and  1.7).

Write $\mathcal{S}^{1}\colon M_{1}\to\wt{N}$ a smooth foliated 
imbedding verifying $\mathfrak{N}_{1} \rondp \mathcal{S}^{1}= \Ide$.
The composition $\mathcal{S} = \mathcal{S}^{1}\rondp 
 \mathcal{S}_{1} \colon M \to \wt{N}$ is a smooth 
foliated imbedding verifying $\mathfrak{N} \rondp \mathcal{S} =
\Ide$ (cf. 1.6).
From Remark 2.3.3 we get that 
$\mathcal{S}_{1}^{*}\kappab_{1}= \kappab$ with $\kappab_{1} \in \left( 
\hiru{H}{1}{M_{1}/\F_{1}}
\right)^{\Z_{2}}$. The above Lemma gives 
$\kappab_{1} =0$. By induction hypothesis  ($\depth \SHhat < \depth 
\SH$) we get that
$\HH_{1}$ is taut. So, the restriction of $\HH_{1}$ to 
$\mathfrak{L}^{-1}(M) = M \times \{-1,1\}$ is also taut. We conclude 
that $\F$ is also taut (cf. 1.6).
\qed

\prgg {\bf Remark}. The proof of the above Theorem shows that the tautness of 
$\F$ and $\wt\HH$ are closely related. In fact,
\begin{eqnarray*}
\hbox{The foliation $\F$ is taut} & \Longleftrightarrow & \hbox{The foliation 
$\wt{\HH}$ is taut.}\\
\Updownarrow \hspace{2cm}& & \hspace{2cm}\Updownarrow  \\
\hbox{The tautness class $\kappab \in \hiru{H}{1}{M/\F}$ vanishes} 
& \Longleftrightarrow & \hbox{
The tautness class $\kappabt \in \hiru{H}{1}{\wt{N}/\wt{\HH}}$ 
vanishes.}
\end{eqnarray*}

\prg{\bf Second characterization of tautness: the bottom group
$\lau{H}{0}{\kappa}{\mf}$.}

We  
give a characterization of the tautness of $\F$ using
$\lau{H}{0}{\kappa}{\mf}$.  Notice that, in the compact case, this result
comes directly from \refp{Masa} and the PoincarŽ Duality of
\cite{KT,KT2} when $M$ is oriented and $\F$ is transversally oriented. 
In fact, we shall not need these orientability conditions.

\bT
\label{T5}
Let $M$ be a manifold endowed with a CERF $\F$. Consider $\mu$ a 
D-metric on $M$.
Then, the following two statements are equivalent:

\Zati The foliation $\F$ is taut.

\zati The cohomology group $\lau{H}{0}{\kappa_{\mu}}{M/\F}$ is $\R $ 
(cf. 2.1).
\bigskip

\nt Otherwise, $\lau{H}{0}{\kappa_{\mu}}{\mf}= 0$.
\eT
\pro
We proceed in two steps.

\medskip

$(a) \Rightarrow (b)$. If $\F$ is taut then $\kappab = 
[\kappa_{\mu}]=0$. So, 
$\lau{H}{0}{\kappa_{\mu}}{\mf} \cong \hiru{H}{0}{\mf} = \R$.

\smallskip

$(b) \Rightarrow (a)$.  If $\lau{H}{0}{\kappa_{\mu}}{\mf}  \neq 0$ 
then there exists a function $0 \neq f \in \hiru{\Om}{0}{\mf}$ with $df = f 
\kappa_{\mu}$. The set $Z(f) = f^{-1}(0)$  is clearly a closed subset of 
$M$. Let us see that it is also an open subset. Take $x \in Z(f)$ and 
consider a contractible open subset $U \subset M$ containing $x$. So, 
there exists a smooth map $g \colon U \to \R$ with $\kappa_{\mu} = 
dg $ on $U$. The calculation
$$
d(f e^{-g}) = e^{-g} df - fe^{-g} dg = e^{-g} f \kappa_{\mu} - e^{-g} f 
\kappa_{\mu} =0
$$
shows that $f e^{-g}$ is constant on $U$. Since $x \in Z(f)$ then $f 
\equiv 0$ on $U$ and therefore $x \in U \subset Z(f)$. We get that 
$Z(f)$ is an open subset.

By connectedness we have that $Z(f) = 
\emptyset$ and $|f|$ is a smooth function. From the equality 
${\displaystyle
d(\log |f|) = \frac{1}{f} df = \kappa_{\mu}
}$
we conclude that $\kappab =0$ and then $\F$ is taut.

\medskip

Notice that we have also proved:  $\lau{H}{0}{\kappa_{\mu}}{\mf}  \neq 0
\Rightarrow \lau{H}{0}{\kappa_{\mu}}{\mf}  = \R$. This ends the proof.
\qed

\prg {\bf Third characterization of tautness: the top group
$\lau{H}{n}{c}{\mf}$.}

We  
give a characterization of the tautness of $\F$ by using  
$\lau{H}{n}{c}{\mf}$. We lift the question to a 
Molino's desingularisation of $\F$, where the result is known.
But we need to formulate an orientability condition on $\F$. 

\bL
\label{TO}
$$
\hbox{ The foliation $\F$ is transversally orientable } \Longleftrightarrow 
\hbox{ The foliation $\wt\HH$ is transversally orientable. }
$$
\eL
\pro
Since $\mathcal{S} \colon 
M\to M_{1}$ is a smooth foliated imbedding 
then we get ``$\Leftarrow$''.

Consider $\mathfrak{O}$ a transverse orientation on $(M,\F)$.  The
tubular neighborhood $(T_0 \cap M,\F)$ inherits the transverse
orientation $\mathfrak{O}$.  Since $((D_0\cap M)\times ]0,1[,\F \times
\mathcal{I})$ is foliated diffeomorphic to $(T_0 \cap M,\F)$ then it
inherits a transverse orientation, written $\mathfrak{O}$.  This
transverse orientation induces on the product $((D_0\cap M)\times
]-1,1[,\F \times \mathcal{I})$ a transverse orientation, written
$\mathfrak{O}$.  Notice that the involution $(x,t) \mapsto (x,-t)$
reverses the orientation $\mathfrak{O}$.

Since 
$$
M_{1} = 
\left\{
\Big( (D_0\cap M)\times ]-1,1[\Big) \coprod \Big( M \times \{ 
-1,1\}\Big)
\right\} \Big/ \sim ,
$$
where $(z,t) \sim ( 2|t|\cdot z , t/|t|)$,
then it suffices to define on $M_1$ the transverse orientation
$\mathfrak{O}_{1}$ by:
\begin{itemize}
    \item[*] \ \ $\mathfrak{O}$ on $(D_0\cap M)\times 
    ]-1,1[$,
    \item[*] \ \ $\mathfrak{O}$ on $M \times \{ 1 \}$, and
    \item[*] $-\mathfrak{O}$ on $M \times \{ - 1 \}$.
    \end{itemize}
This gives ``$\Rightarrow$''. \qed

Before passing to the third characterization, we need two computational
Lemmas. 

\bL
\label{E1}
The inclusion $$\lau{\Om}{*}{c}{((D_0 \cap M) \times ]0,1[)/\F \times 
\mathcal{I}} \hookrightarrow \lau{\Om}{*}{c}{((D_0 \cap M)  \times 
]-1,1[)/\F \times 
\mathcal{I}}$$
induces an  isomorphism in cohomology.
\eL
 \pro
For the sake of simplicity, we write $E = (D_0 \cap M)$. 
Let $f \in \lau{\Om}{0}{c}{]-1,1[}$ be a function with
 $f\equiv 0$ on $]-1,1/3]$ and $ f \equiv 1$ on $[2/3,1[$. So, $df \in \lau{\Om}{1}{c}{]0,1[} \subset 
 \lau{\Om}{1}{c}{]-1,1[}$. 
Lemma will be proved if we show that the assignment $[\gamma] \mapsto
[df \wedge \gamma ]$ establishes the isomorphisms of degree $+1$
 $$
 \lau{H}{*}{c}{E/\F} \cong \lau{H}{*}{c}{(E \times ]0,1[)/\F\times 
 \mathcal{I}} 
 \ \ \hbox{ and } \ \ 
 \lau{H}{*}{c}{E/\F} \cong  \lau{H}{*}{c}{(E \times ]-1,1[)/\F\times 
 \mathcal{I}}.
 $$

Let us prove the first one (the second one is proved in the same way). 
Consider the following differential complexes:
\begin{itemize}
 \item[-] $\hiru{A}{*}{ ]0,3/4[}  = 
     \left\{ \om \in \hiru{\Om}{*}{E \times  ]0,3/4[/ \F\times 
 \mathcal{I}}  \Big/ \ \left[
 \begin{array}{l}
 \supp \om \subset K \times [c,3/4[\\
 \hbox{for a compact } K \subset E 
 \hbox{ and } 0 <c < 3/4
 \end{array} 
 \right.
 \right\}.
 $
 \item[-] $\hiru{A}{*}{ ]1/4,1[}  = 
 \left\{ \om \in \hiru{\Om}{*}{E \times  ]1/4,1[/ \F\times 
 \mathcal{I}}  \Big/ \ \left[
 \begin{array}{l}
 \supp \om \subset K \times ]1/4,c]\\
 \hbox{for a compact } K \subset E 
 \hbox{ and } 1/4 < c <1
 \end{array} 
 \right.
 \right\}.
 $
 \item[-] $\hiru{A}{*}{ ]1/4,3/4[}  = 
 \left\{ \om \in \lau{\Om}{*}{}{E \times ]1/4,3/4[/ \F\times 
 \mathcal{I}} \Big/ \ 
 \left[\begin{array}{l}
 \supp \om \subset K \times ]1/4,3/4[\\
 \hbox{for a compact } K \subset E
 \end{array}
 \right.
 \right\}.
 $
 \end{itemize}
Proceeding as in \refp{mv} we get the short exact sequence
 $$
 0 \TO \lau{\Om}{*}{c}{E \times ]0,1[/\F \times \mathcal{I}}
 \TO
 \hiru{A}{*}{ ]0,3/4[}  \oplus \hiru{A}{*}{ ]1/4,1[} 
 \TO
 \hiru{A}{*}{ ]1/4,3/4[}
 \TO 0.
 $$
 The associated long exact sequence is
\begin{eqnarray*}
\cdots &\to &
\hiru{H}{i-1}{\hiru{A}{*}{ ]1/4,3/4[} } 
\stackrel{\delta}{\to}
\lau{H}{i}{c}{E \times ]0,1[/\F \times \mathcal{I}}
 \to 
 \hiru{H}{i}{\hiru{A}{*}{ ]0,3/4[} } \oplus \\
 &\oplus& \hiru{H}{i}{\hiru{A}{*}{ 
 ]1/4,1[} }
\to
 \hiru{H}{I}{\hiru{A}{*}{ ]1/4,3/4[}}
\to \cdots,
\end{eqnarray*}
 where the connecting morphism is
 $
 \delta([\om]) = [df \wedge \om].
 $
\medskip

Before executing the calculation let us introduce some notation. 
       Let $\beta$ be a differential form on $\hiru{\Om}{i}{D_S \times 
       ]a,b[}$
       which does not 
       include the $dt$ factor. 
By $\Int{-}{c}\beta (s)\wedge ds$ and
       $\Int{c}{-}\beta  (s)\wedge ds$ we denote the forms on $\hiru{\Om}{i}{D_S \times 
       ]a,b[}$ 
       obtained from $\beta$ 
       by integration with respect to $s$, that is,
       $\left(\Int{-}{c}\beta (s)\wedge ds \right) 
       (x,t)(\vec{v}_{1}, 
       \ldots ,\vec{v}_{i}) = \Int{t}{c}(\beta (x,s)(\vec{v}_{1}, 
       \ldots ,\vec{v}_{i}) ) \ ds $
      and  on the other hand $\left( \Int{c}{-}\beta (s)\wedge ds 
       \right)(x,t)(\vec{v}_{1}, 
       \ldots ,\vec{v}_{i}) )= \Int{c}{t}(\beta 
       (x,s)(\vec{v}_{1}, 
       \ldots ,\vec{v}_{i}) ) \,  ds$
       where  $c \in ]a,b[$, $(x,t) \in D_S \times ]a,b[$ and 
       $(\vec{v}_{1}, 
       \ldots ,\vec{v}_{i}) \in T_{(x,t)}( D_{S} \times ]a,b[) $ .

\medskip  

\underline{\em (i) Computing $\delta$}. 

\smallskip

 Each differential form $\om \in \hiru{A}{*}{Interval}$ 
 can be written $\om  = \alpha + \beta \wedge dt$ where $\alpha$ and 
 $\beta$ do not 
 contain $dt$. 

 Consider  a cycle $\om  = \alpha + \beta \wedge dt \in \hiru{A}{i}{ 
 ]0,3/4[}$ with $\supp \om \subset K \times [c,3/4[$ for a 
 compact 
 $K \subset E$ and $0 < c < 3/4$.
 We have 
 $
 \om = \alpha(c/2)   - d \left( \Int{-}{c/2} \beta (s) \wedge  ds\right) = 
 -d \left( \Int{-}{c/2} \beta (s) \wedge  ds\right).
 $
 Since $ \supp \Int{-}{c/2} \beta (s) 
 \wedge ds \subset K \times [c,3/4[$ then we get
 $
 \hiru{H}{*}{\hiru{A}{\cdot}{ ]0,3/4[} } =0$. In the same way, we 
 get  $
 \hiru{H}{*}{\hiru{A}{\cdot}{ ]1/4,1[} } =0$. 

 \smallskip
 
 We conclude that $\delta$ is an isomorphism.

 \medskip  

 \underline{\em (i) Computing $\lau{H}{*}{c}{(E \times ]0,1[)/\F
 \times \mathcal{I}}$}. 

 \smallskip

 Consider  a cycle $\om  = \alpha + \beta \wedge dt \in \hiru{A}{*}{ 
]1/4,3/4[}$. 
 We have 
 $
 \om = \alpha(1/2)   +d \left( \Int{1/2}{-} \beta (s) \wedge  ds\right). 
 $
 Notice that $ \supp \Int{1/2}{-}  \beta (s) 
 \wedge ds \subset K \times ]1/4,3/4[$ and $\supp \alpha (1/2) \subset K$.
A standard procedure shows that the operator $\Delta \colon 
\hiru{H}{*}{\hiru{A}{\cdot}{ ]1/4,3/4[} }
\to \lau{H}{*-1}{c}{E/\mathcal{F}}$, defined by $\Delta([\om]) = 
[\alpha (1/2)]$, is  an isomorphism.
The inverse is 
$\Delta^{-1}([\gamma]) = [\gamma]$.

So, the composition $\delta \rondp\Delta^{-1}  \colon \lau{H}{*-1}{c}{E/\mathcal{F}}
\to \lau{H}{*}{c}{(E \times ]0,1[)/\F
\times \mathcal{I}}$ is an isomorphism. It is exactly the operator:
$[\gamma] \mapsto [df\wedge \gamma]$.
 \qed

The reason why we use the $(-\Z_{2})$-invariant classes in the next Lemma 
instead of
the more natural $\Z_{2}$-invariant classes is 
the following: we have 
$
 \left(\lau{H}{*}{c}{]\! - \! 1,0[ \cup ]0,1[/\mathcal{I}}\right)^{-\Z_{2}} \cong 
 \left( \lau{H}{*}{c}{]\! - \! 1,1[)/\mathcal{I}}\right)^{-\Z_{2}}
 $
but also
$
 \left(\lau{H}{*}{c}{]-1,0[ \cup ]0,1[/\mathcal{I}}\right)^{\Z_{2}} 
 \not\cong 
 \left( \lau{H}{*}{c}{]-1,1[)/\mathcal{I}}\right)^{\Z_{2}}.
 $
  \bL
  \label{E2}
  The inclusion $$\lau{\Om}{*}{c}{((D_0 \cap M) \times (]-1,0[ \cup ]0,1[))/\F\times 
  \mathcal{I}} \hookrightarrow \lau{\Om}{*}{c}{((D_0 \cap M) \times ]-1,1[)/\F \times 
  \mathcal{I}}$$
  induces the isomorphism
  $$
  \left(\lau{H}{*}{c}{((D_0 \cap M) \times (]-1,0[ \cup ]0,1[))/\F\times 
  \mathcal{I}}\right)^{-\Z_{2}} \cong 
  \left( \lau{H}{*}{c}{((D_0 \cap M) \times ]-1,1[)/\F\times 
  \mathcal{I}}\right)^{-\Z_{2}}.
  $$
  \eL
  \pro
  For the sake of simplicity, we write $E = (D_0 \cap M)$. 
  Consider $f \colon ]-1,1[ \to \R$ a smooth function with 
   $f\equiv 0$ on $]-1,1/3]$ and $ f \equiv 1$ on $[2/3,1[$. 
The function $f \rondp \sigma + f -1 \colon ]-1,1[ \to \R$ is a 
smooth function whose support is in $[-2/3,2/3]$. So, 
$[df] = -[d(f\rondp \sigma)] \in \lau{H}{*}{c}{]-1,1[}$.
The above Lemma gives
  \begin{itemize}
      \item[-] $\lau{H}{*}{c}{E/\mathcal{F}} \cong \lau{H}{*}{c}{(E \times ]-1,1[)/\F\times 
\mathcal{I}}$ by $[\gamma ] \mapsto [df] \wedge [\gamma]$,
      \item[-]$\lau{H}{*}{c}{E/\mathcal{F}} 
      \oplus \lau{H}{*}{c}{E/\mathcal{F}}
      \cong \lau{H}{*}{c}{(E \times 
      (]-1,0[ \cup ]0,1[))/\F\times 
\mathcal{I}}
$ by $([\gamma _{1}],[\gamma_{2}]) \mapsto 
[d(f \rondp \sigma)] \wedge [\gamma_{1}] +
[df] \wedge \gamma_{2}$.
      \end{itemize}
 Notice that 
 $
 \zeta \cdot df = - \sigma^{*} df = - d(f \rondp \sigma)$ on $]-1,1[$ 
 and $]-1,0[ \cup ]0,1[$.
 This gives the isomorphisms of degree $+1$:
 \begin{itemize}
 \item[-] $\lau{H}{*}{c}{E/\mathcal{F}} \cong 
 \left(\lau{H}{*}{c}{(E \times ]-1,1[)/\F\times 
\mathcal{I}}\right)^{-\Z_{2}}$ by $[\gamma ] \mapsto [df] \wedge [\gamma]$,

   \item[-]$\lau{H}{*}{c}{E/\mathcal{F}} 
   \cong \left(\lau{H}{*}{c}{(E \times 
   (]-1,0[ \cup ]0,1[))/\F\times 
\mathcal{I}}\right)^{-\Z_{2}}$ by $[\gamma]\mapsto 
-\frac{1}{2}[d(f \rondp \sigma)] \wedge [\gamma] +
\frac{1}{2}[df] \wedge [\gamma]$.
   \end{itemize}
   We obtain the result, since on $]-1,1[$ we have
   $
   -\frac{1}{2}[d(f \rondp \sigma)]  +
   \frac{1}{2}[df] = \frac{1}{2}[df]  +
   \frac{1}{2}[df] = [df].
   $
  \qed

The way to lift the question  to the Molino's desingularisation is the
following.

\bL
\label{nulo}
$$
\lau{H}{n}{c}{M_1/\F_{1}} \neq 0 \Longleftrightarrow
\lau{H}{n}{c}{M/\F} \neq 0.
$$
\eL
\pro
We prove the two implications.

\smallskip

\fbox{$\Leftarrow$} From the open covering $\left\{ (D_0\cap M) \times ]-1,1[, 
M \times \{ -1,1\}\right\}$ of $M_1$ we obtain the short 
exact sequence
\begin{eqnarray*}
0 
&\to&
\lau{\Om}{*}{c}{(D_0\cap M) \times (]-1,0[ \cup 
]0,1[)/ \F \times \mathcal{I}} 
\to
\\
&\to& \lau{\Om}{*}{c}{(D_0\cap M) \times ]-1,1[ / \F 
\times \mathcal{I}} 
\oplus
\lau{\Om}{*}{c}{(M \times \{-1,1\})/ \F \times \mathcal{I}} \to
\lau{\Om}{*}{c}{M_1 / \F_1 }  \to 0
\end{eqnarray*}
(same argument as in \refp{mv}). The associated long exact sequence is 
$$\cdots 
\to
\lau{H}{j}{c}{(D_0\cap M) \times (]-1,0[ \cup 
]0,1[)/ \F \times \mathcal{I}} 
\stackrel{(I,J)}{\TO}
\lau{H}{j}{c}{(D_0\cap M) \times ]-1,1[ / \F 
\times \mathcal{I}} \oplus
$$
$$
\lau{H}{j}{c}{(M \times \{-1,1\})/ \F \times \mathcal{I}} \to 
\lau{H}{j}{c}{M_1 / \F_1 }  \to \lau{H}{j+1}{c}{(D_0\cap M) \times (]-1,0[ \cup 
]0,1[)/ \F \times \mathcal{I}} 
\to \cdots
$$

Lemma \ref{E1} gives that $I$ is an onto map. We also have
$$
\lau{H}{n+1}{c}{(D_0\cap M) \times (]-1,0[ \cup 
]0,1[)/ \F \times \mathcal{I}}
=0
$$
for degree reasons.  Therefore
$$
\lau{H}{n}{c}{M  / \F}
\oplus
\lau{H}{n}{c}{M  / \F}
=
\lau{H}{n}{c}{(M  \times \{ -1,1\})/ \F  \times 
\mathcal{I}}  
\to 
\lau{H}{n}{c}{M_1 / \F_1 }  
$$
is a surjective map.  This gives the result.

\smallskip

\fbox{$\Rightarrow$} The above short 
exact sequence is $\Z_{2}$-equivariant. So, we get the long exact sequence
\begin{eqnarray*}
\! \!\cdots \! \!
&\to&
\left(\lau{H}{j}{c}{(D_0\cap M) \times (]-1,0[ \cup 
]0,1[)/ \F \times \mathcal{I}} \right)^{-\Z_{2}}
\stackrel{(I,J)}{\TO}
\\
&\to& \left(\lau{H}{j}{c}{(D_0\cap M) \times ]-1,1[ / \F 
\times \mathcal{I}} \right)^{-\Z_{2}}
\oplus
\left(\lau{H}{j}{c}{(M \times \{-1,1\})/ \F \times \mathcal{I}}\right)^{-\Z_{2}} 
\to \\
&\to&
\left(\lau{H}{j}{c}{M_1 / \F_1 }\right)^{-\Z_{2}} \to \left(\lau{H}{j+1}{c}{(D_0\cap M) \times (]-1,0[ \cup 
]0,1[)/ \F \times \mathcal{I}}\right)^{-\Z_{2}}
\to \cdots
\end{eqnarray*}
Lemma \ref{E2} ensures that $I$ is an isomorphism.  We conclude that
$$
\left(\lau{H}{*}{c}{M_1 / \F_1 }  \right)^{-\Z_{2}}
\! \cong \! 
\left( \lau{H}{*}{c}{(M  \times \{-1,1\})/ \F \times 
\mathcal{I}}\right)^{-\Z_{2}}
\! \cong \! 
\lau{H}{*}{c}{M/ \F } \! \otimes \!  \left( \hiru{H}{0}{-1,1}\right)^{-\Z_{2}}
\! \cong \! 
\lau{H}{*}{c}{M / \F } ,
$$
which ends the proof.\qed

The main result of this section is the following.
\bT
\label{T3}
Let $M$ be a manifold endowed with a CERF $\F$. 
Let us suppose that $\F$ is transversally orientable.
Put $n = \codim \F$.
Then, the two following statements are equivalent:

\Zati The foliation $\F$ is taut.

\zati The cohomology group $\lau{H}{n}{c}{M/\F}$ is $\R $.

\bigskip

\nt Otherwise, $\lau{H}{n}{c}{\mf}= 0$.
\eT
\pro
We proceed by induction on $\depth \SH$. When $\depth \SH =0$ then we 
have a regular foliation and the result comes from  \cite[Theorem 
6.4]{Suso}. For the induction step, we can suppose that the result is 
true for $\SHhat$, since $\depth \SHhat < \depth \SH$. We have in 
particular $\lau{H}{n}{c}{M_{1}/\F_{1}} = \R$ or $0$. So, Lemma  
\ref{nulo} 
ensures that
$$
\lau{H}{n}{c}{M_{1}/\F_{1}} = \R \ (\hbox{resp. } 0)
\Longleftrightarrow
\lau{H}{n}{c}{M/\F} = \R \ (\hbox{resp. } 0).
$$

By Lemma \ref{TO} the
foliation $\wt\HH$  is trans\-versally 
orientable. Since $(\wt{H},\wt{\HH})$ is a Molino's desingularisation 
of $(N_{1},\HH_{1})$ then the foliation $\HH_{1}$ is also trans\-versally 
orientable.  So, we have
$$
\F \hbox{ is taut } \stackrel{2.4.3}{\Longleftrightarrow}
\wt{\HH} \hbox{ is taut } \stackrel{2.4.3}{\Longleftrightarrow}
\F_{1} \hbox{ is taut } \stackrel{Induc}{\Longleftrightarrow}
\lau{H}{n}{c}{M_1/\F_1} = \R \Longleftrightarrow
\lau{H}{n}{c}{M/\F} = \R ,
$$
which ends the proof.\qed

\prg {\bf Reading the tautness of $\F$ on $N$ and $\wt{N}$}. The tautness of 
$\F$ and $\wt\HH$ are closely related. We have already seen
$$
\lau{H}{n}{c}{\mf} = \R \Longleftrightarrow 
\hiru{H}{n}{\wt{N}/{\wt{\HH}}} = \R.
$$
The situation at the zipper level is more complicated.
In fact, we don't have the expected equivalence\footnote{The basic cohomology of 
an SRF is defined as in 2.1. } $
\lau{H}{n}{c}{\mf} = \R \Leftrightarrow 
\hiru{H}{n}{N/\HH} = \R $ 
(for instance, in Example 1.4 we have  $\lau{H}{n}{c}{\mf} = \R$ and
$\hiru{H}{n}{N/\HH} = 0$).  This situation comes from the existence of
boundary strata but we are going to prove that it is the only
obstruction.

\bL
\label{T}
Consider $(T_S,\tau,S)$ a foliated tubular neighborhood of a 
singular stratum $S \in \SH^{^{sin}}$. If $S$ is not a boundary stratum then
the inclusion $(T_S\menos S )\hookrightarrow  T_S$ induces the isomorphism 
$
\iota \colon\lau{H}{n}{c}{(T_S\menos S)/\mathcal{\HH}} \to
\lau{H}{n}{c}{T_S/\mathcal{\HH}}.
$
\eL
\pro
We proceed in two steps.

\Zati  \underline{\em Approaching $\lau{H}{n}{c}{T_S/\F}$}. Consider the complexes
\begin{itemize}
    \item[-] $\hiru{A}{*}{ [0,3/4[}  = 
    \left\{ \om \in \lau{\Om}{*}{}{\rho_S^{-1}( [0,3/4[ )/ \mathcal{\HH} } \ \Big/ \ \left[
\begin{array}{l}
\supp \om \subset \tau^{-1}(K) \\
\hbox{for a compact } K \subset S.
\end{array} 
\right.
\right\}.
$
\item[-] $\hiru{A}{*}{ ]1/4,1[}  = 
\left\{ \om \in \lau{\Om}{*}{}{\rho_S^{-1}( ]1/4,1[)/ \mathcal{\HH} } \ \Big/ \ 
\left[\begin{array}{l}
\supp \om \subset\tau^{-1}(K) \cap \rho_S^{-1}( ]1/4,\epsilon])\\
\hbox{ for a compact } K \subset S 
\hbox{ and } 1/4 < \epsilon < 1.
\end{array}\right.
\right\}.
$
\item[-] $\hiru{A}{*}{ ]1/4,3/4[}  = 
\left\{ \om \in \lau{\Om}{*}{}{\rho_S^{-1}( ]1/4,3/4[)/ \mathcal{\HH} } \ \Big/ \ 
\left[\begin{array}{l}
\supp \om \subset\tau^{-1}(K)  \\
\hbox{ for a compact } K \subset S .
\end{array}
\right.
\right\}.
$
\end{itemize}

The short exact sequence
$$
0 \TO \lau{\Om}{*}{c}{T_S/ \mathcal{\HH} } 
\TO
\hiru{A}{*}{ [0,3/4[}  \oplus \hiru{A}{*}{ ]1/4,1[} 
\TO
\hiru{A}{*}{ ]1/4,3/4[}
\TO 0
$$
produces the long exact sequence
\begin{eqnarray*}
\hiru{H}{n-1}{\hiru{A}{*}{ [0,3/4[} } \oplus \hiru{H}{n-1}{\hiru{A}{*}{ ]1/4,1[} }&\TO&
\hiru{H}{n-1}{\hiru{A}{*}{ ]1/4,3/4[}}
\stackrel{\delta}{\TO}
\lau{H}{n}{c}{T_S/\mathcal{\HH}}
\TO
\\
&\TO&
\hiru{H}{n}{\hiru{A}{*}{ [0,3/4[} } \oplus \hiru{H}{n}{\hiru{A}{*}{ ]1/4,1[} }
\end{eqnarray*}
where the connecting morphism is
$
\delta([\om ]Ê) = [df \wedge \om]
$
for a smooth map $f\colon [0,1[ \to [0,1]$ with $f \equiv 0$ on 
$[0,3/8]$ and $f\equiv 1$ on $[5/8,1[$.

\zati \underline{\em Computing $\hiru{H}{*}{\hiru{A}{\cdot}{Interval}}$}. 
Consider  a cycle $\om  = \alpha + \beta \wedge dt \in \hiru{A}{*}{ 
]1/4,1[}$, where $\alpha$ and 
$\beta$ do not 
contain $dt$. 
We have 
$
\om = 
-d \left( \Int{-}{1} \beta (s) \wedge  ds\right).
$
Since $ \supp \Int{-}{1} \beta (s) 
\wedge ds \subset \tau^{-1}(K) \times ]1/4,\epsilon]$ we get
$
\hiru{H}{*}{\hiru{A}{\cdot}{ ]1/4,1[} } =0$.

The contraction 
$H_{S} \colon \rho_{S}^{-1}([0,3/4[) \times [0,1] \to \rho_{S}^{-1}([0,3/4[) $ (cf. 1.4) 
is a foliated proper homotopy between $\rho_{S}^{-1}([0,3/4[) $ and 
$S$. This homotopy preserves the fibers of $\tau$. Then 
$
\hiru{H}{*}{\hiru{A}{\cdot}{ [0,3/4[)} }
\cong
\lau{H}{*}{c }{S/\HH_{S}}.
$
Since the stratum $S$ is a not boundary stratum, then 
$\codim_{S}\HH_{S} < n-1$ (cf. 1.4). Then, by degree reasons, we get 
$
\hiru{H}{n-1}{\hiru{A}{*}{ [0,3/4[)} } =\hiru{H}{n}{\hiru{A}{*}{ [0,3/4[)} }=0$.

Since  $ \rho_{S}^{-1}(]1/4,3/4[) = D_{S}\times ]1/4,3/4[$, then we 
have that the assignment $[\gamma] \mapsto [\gamma \wedge df]$ induces the 
isomorphism 
$
\lau{H}{*}{c }{D_{S}/\HH}
\cong
\hiru{H}{*}{\hiru{A}{\cdot}{ ]1/4,3/4[} }
$
(cf. proof of the Lemma \ref{E1}).

\zati {\em Last step.} The above exact sequence gives that
$
\delta \colon \lau{H}{n-1}{c}{D_{S}/\mathcal{\HH}}
\to
\lau{H}{n}{c}{T_{S}/\mathcal{\HH}}
$,
defined by $\delta ([\om]) = [df \wedge \om]$, is an isomorphism.
We have completed the proof as $\delta \colon
\lau{H}{n-1}{c}{D_{S}/\HH} \to \lau{H}{n}{c}{(T_{S}\menos S/\HH)
\equiv (D_{S} \times ]0,1[)/ \mathcal{H} \times \mathcal{I}}$ is an
isomorphism (cf.  proof of the Lemma \ref{E1}).  \qed

The tautness of $\F$ can be read on a zipper as follows.
\bT
\label{T4}
Let $M$ be a manifold endowed with a CERF $\F$, which is 
transversally oriented.  Consider $(N,\HH)$ a 
zipper of $\F$. Let us suppose that $\SH$ does not possess any
boundary stratum.
Put $n = \codim \F$.
Then, the following two statements are equivalent:

\Zati The foliation $\F$ is taut.

\zati The cohomology group $\hiru{H}{n}{N/\HH}$ is $\R $.

\bigskip

\nt Otherwise, $\hiru{H}{n}{N/\HH} = 0$.
\eT
\pro
It suffices to prove that $ \hiru{H}{n}{N/\HH} = \lau{H}{n}{c}{\mf}$ (cf. 
Theorem \ref{T3}).  Consider $\Sigma_i $ as in \ref{reduc}.  We have
$N \menos \Sigma_{-1} = N$ and $N \menos \Sigma_{m-1} = M$, where $m =
\dim M$.  We will get the result if we prove that $\lau{H}{n}{c}{(N\menos
\Sigma_{i})/\HH} \cong \lau{H}{n}{c}{(N\menos \Sigma_{i-1})/\HH}$, for
$i \in \{ 0, \ldots , m-1\}$.

From the open covering
$
{\displaystyle
\left\{ N\menos\Sigma_i, \bigcup_{\dim S =i}T_{S}\right\}
}
$
of $N\menos\Sigma_{i-1}$, 
we obtain the Mayer-Vietoris sequence
$$
\bigoplus_{\dim S =i} \lau{H}{n}{c}{(T_{S}\menos S)/\HH} \TO
\lau{H}{n}{c}{(N\menos\Sigma_i)/\HH}
\oplus
\bigoplus_{\dim S =i} \lau{H}{n}{c}{T_{S}/\HH} 
\TO
\lau{H}{n}{c}{(N\menos\Sigma_{i-1})/\HH} \to 0.
$$
Now, the Lemma \ref{T} gives the result.
\qed


When boundary strata appear, it can be shown that this Theorem 
still holds considering $\lau{H}{n}{}{N/\HH, \partial(N/\HH)}$ in (b). 
Here, $\partial(N/\HH) = \cup \{ S \ | \ S \preceq S' \hbox{ for a 
boundary stratum } S'\}$ and the relative basic cohomology must be 
understood in a suitable way.

\section{Appendix}
We prove the existence of foliated Thom-Mather systems for an SRF $\HH$
defined on a compact manifold $N$ announced in Proposition \ref{TM}. 
First of all, we need a more accurate presentation of the Molino's
desingularisation of 1.7.  We fix an adapted metric $\mu$ on $N$.

\prg {\bf Molino's desingularisation}. A {\em Molino's desingularisation} 
of $(N,\HH,\mu)$ (see \cite{Mo0}) is a sequence $\mathfrak{N}$ of blow
ups (cf.  1.6):
\begin{equation}
    \label{bup}
(N_{p},\HH_{p},\mu_{p}) \stackrel{\mathcal{L}_{p}}{\TO}
(N_{p-1},\HH_{p-1},\mu_{p-1})  \TO \cdots \TO
(N_{1},\HH_{1},\mu_{1}) \stackrel{\mathcal{L}_{1}}{\TO}
(N_{0},\HH_{0},\mu_{0}) 
\end{equation}
with
$(N_{0},\HH_{0},\mu_{0})   = (N,\F,\mu)$ and $ (N_{i},\HH_{i},\mu_{i}) \equiv 
    (\widehat{N_{i-1}},\widehat{\HH_{i-1}},\widehat{\mu_{i-1}})$, 
    for $i \in \{ 1, \ldots , p\}$\footnote{The map $\xi$  necessary 
    for the construction of the blow up (cf. 1.6)    is supposed to be the same 
    for each $i$.}. 
    Here, $p = \depth \SF$. The triple
    $(\wt{N},\wt{\HH},\wt{\mu}) \equiv 
    (N_{p},\HH_{p},\mu_{p}) $
    is a regular riemannian foliated manifold.
    
    Notice that the blow up of 
    $(N \times \R,\F \times \I, \mu + dt^{2})$ is just
$$
    (N_{p} \times \R,\HH_{p} \times \I,\mu_{p} + dt^{2}) \stackrel{\mathcal{L}_{p}}{\TO}
    (N_{p-1} \times \R,\HH_{p-1}\times \I,\mu_{p-1}+ dt^{2})  \TO \cdots 
    \stackrel{\mathcal{L}_{1}}{\TO}
    (N_{0} \times \R,\HH_{0}\times \I,\mu_{0}+ dt^{2})
$$
(cf. 1.6.1 (a)).

\prg {\bf Construction of the foliated Thom-Mather system.} We first
construct a tubular neighborhood of each singular stratum and then, we
prove that these neighborhoods are satisfy the compatibility condition
\refp{TM1}.

\prgg Take $S \in \SH$ a singular stratum and let $i \in \N$  be its depth.   
 Take the partial desingularisation $\mathfrak{N}_{i}
= \mathfrak{L}_{1} \rondp \mathfrak{L}_{2} \rondp \cdots \rondp
\mathfrak{L}_{i} \colon (N_{i},\HH_{i},\mu_{i}) \to (N,\F,\mu)$. 
Notice that each restriction
\begin{equation}
     \label{trivial}
 \mathfrak{N}_{i}\colon \mathfrak{N}_{i}^{-1}(N\menos \Sigma_{i-1}) 
 \to N\menos \Sigma_{i-1}
 \end{equation}
 is a trivial foliated smooth bundle (cf. 1.6.1 (c),(e)). We fix 
 $\mathfrak{s}_{i}\colon N\menos \Sigma_{i-1} \to 
 \mathfrak{N}_{i}^{-1}(N\menos \Sigma_{i-1})$ a smooth foliated section 
 of \refp{trivial}.  For $i=0$ we put $\mathfrak{N}_{i} =
 \mathfrak{s}_{i} = \Ide \colon N \to N$.

The stratum $S$ is a proper submanifold of the foliated riemannian manifold
$(N\menos\Sigma_{i-1},\HH_{i},\mathfrak{s}_{i}^{*}\mu_{i})$.  Let 
$(T_{S}, \tau_{S},S)$ be  a foliated tubular neighborhood of $S$
constructed as in 1.4.  It is also a foliated tubular neighborhood of
$S$ in $(N,\HH)$.  We need to shrink it in order to assure the
Thom-Mather compatibility of these neighborhoods.

Following 1.6.1 (e),
there exists a family of strata $\{ S^j \in \SHj \ | \ j \in \{ 0,
\ldots i\}\}$ such that:
\begin{itemize}
    \item[+] $S^{0} = S$, 
    \item[+] $(S^j)^{\ \widehat{}}  = S^{j+1}$ for $j \in \{0, 
    \ldots , i-1\}$ and
    \item[+] $S^i$ is a minimal stratum of $\SHi$.
    \end{itemize}
Notice that $\mathfrak{s}_{i}(S) \subset S^i$. By construction, there exists a
foliated tubular neighborhood $(U_{S^i},\nu_{S^i},S^i)$ of $S^i$ on
$(N_i,\HH_{i})$ such that its restriction to $S^i$ is isomorphic to
$(T_{S}, \tau_{S},S)$ through $\mathfrak{s}_{i}$.
We can take these neighborhoods small enough such that:
\begin{equation}
    \left( \mathcal{L}_{i+1}^{-1}(U_{S^i}), \HH_{i+1},\mu_{i+1}\right) =
    \left( D_{S^i} \times ]-1,1[, \HH_{i+1}
    \times \mathcal{I},\mu_{i}|_{D_{S^i} }+ dt^{2}\right), \end{equation}
    where $D_{S^i} $ is the core of $(U_{S^i},\nu_{S^i},S^i)$ (cf.  1.6).

    \prgg We prove now that the family $\{T_{S} \ | \ S \in \SF^{^{sin}}\}$ which
    we have constructed is a foliated Thom-Mather
    system of $(N,\HH)$.  Now, we consider two singular strata $S,S' \in
    \SF^{^{sin}}$ with $S \prec S'$ and we prove the property \refp{TM1},
    which will lead to \refp{TM2}.
    
    Without loss of generality we can suppose that $S$ is a
    minimal stratum.  The Thom-Mather conditions involve just $T_{S}$ and
    not the whole manifold $N$.  So, we can suppose $N=T_{S}$.  The blow up $
    \mathfrak{L} \colon (\wh{N},\wh{\HH},\wh{\mu}) \to (N,\HH,\mu) $
    becomes the map $ \mathfrak{L} \colon (D_{S} \times ]-1,1[ , \HH
    \times \mathcal{I},\mu|_{D_{S}} + dt^{2})\TO (T_{S},\HH,\mu), $
    defined by $\mathfrak{L}(u,t) = 2|t| \cdot u$.  Recall that the
    restriction
$$
\mathfrak{L} \colon (D_{S} \times (]-1,1[\menos \{ 0 \}) , \HH \times
\mathcal{I})\TO ((T_{S}\menos S),\HH)$$ is a foliated diffeomorphism.
We check the property \refp{TM1} at the $(D_{S} \times
(]-1,1[\menos \{ 0 \}))$-level.  This makes sense since $S \cap
T_{S}\cap T_{S'} = \emptyset$.

The foliated tubular neighborhood $(T_{S'},\tau_{S'},S')$ can be
described as follows.  There exists a foliated tubular neighborhood
$(T_{S''},\tau_{S''},S'')$ of the stratum $S'' = D_{S}\cap S' \in \SHDS$
such that 
$(T_{S'},\tau_{S'},S')$ can be identified with $ (T_{S''} \times
(]-1,1[\menos \{ 0 \}),\tau_{S''} \times \Ide,S''\times (]-1,1[\menos \{ 0 \}))$ through $
\mathfrak{L}$
(cf.  3.1).  By using the foliated diffeomorphism $
\mathfrak{L}$, the condition \refp{TM1} is clear since:
\begin{itemize}
      \item[+] $T_{S}\cap S'$ becomes $S'' \times (]-1,1[\menos \{ 0 \})$,
    \item[+]  $T_{S}\cap T_{S'}$ becomes $T_{S''} \times (]-1,1[\menos
    \{ 0 \})$, 
    \item[+] $\tau_{S'}$ becomes the map $(z,t) \mapsto (\tau_{S''}(z),t)$, and
    \item[+] $\rho_{S}$ becomes the map $(z,t) \mapsto |t|$.
    \qed
    \end{itemize}

\end{document}